\documentclass[a4paper,10pt]{article}
\usepackage{amssymb,amsmath}%,showlabels}

\newcommand{\x}{\b{x}}
\renewcommand{\b}[1]{\mathbf{#1}}
\newcommand{\bQ}{\b{Q}}
\newcommand{\cl}[1]{\mathcal{#1}}
\newcommand{\cQ}{\cl{Q}}
\newcommand{\cR}{\cl{R}}

\newcommand{\Proj}{\mathbb{P}}
\newcommand{\M}{(Q_3L)^{\perp}}

\newcommand{\Q}{\mathbb{Q}}

\newcommand{\ep}{\varepsilon}
\newcommand{\beql}[1]{\begin{equation}\label{#1}}
\newcommand{\eeq}{\end{equation}}

\newtheorem{theorem}{Theorem}[section]
\newtheorem{definition}{Definition}
\numberwithin{equation}{section}
\newtheorem{lemma}{Lemma}[section]

\begin{document}
\title{Rational Points on the Intersection of Three Quadrics}
\author{D.R. Heath-Brown\\Mathematical Institute, Oxford}
\date{}
\maketitle
\section{Introduction}

Let $Q_1(\x),Q_2(\x),Q_3(\x)\in K[\x]$ be three quadratic forms in $n$
variables $\x=(x_1,\ldots,x_n)$, defined over a number field $K$. This paper
will be concerned with the Hasse principle and weak approximation for points
over $K$ on the intersection 
\[\cR:\; Q_1=Q_2=Q_3=0.\]

We begin by reviewing the situation for individual quadrics and pairs
of quadrics.
In the case of a single quadratic form the Hasse principle is always valid,
while weak approximation holds for all nonsingular 
quadratic forms in $n\ge 3$
variables. When one has a pair of forms the Hasse principle may fail, even
when the variety $Q_1=Q_2=0$ is nonsingular, as is shown by the example 
\[Q_1=X_1X_2-(X_3^2-5X_4^2),\;\;\;
Q_2=(X_1+X_2)(X_1+2X_2)-(X_3^2-5X_5^2)\]
(with $K=\Q$) due to Birch and Swinnerton-Dyer \cite{BSD}.
However it is known that for a nonsingular intersection defined by a pair
of forms in 8 or more variables both the Hasse principle and weak
approximation hold, see Heath-Brown \cite[Theorem 1]{HB8}. One cannot
dispense with the smoothness condition here.  The example
\[6X_1^2-X_2^2-X_3^2=X_4^2+\ldots+X_n^2=0 \]
over the field
$\Q$, has points in every completion, as soon as $n\ge 6$, but has no
rational points. 

The difficulty with the above example arises from the real completions, and
it was shown by 
Colliot-Th\'{e}l\`{e}ne, Sansuc and Swinnerton-Dyer \cite[Theorem C]{CTSSD}
that any pair of quadratic forms defined over a totally complex number field
$K$ has a common zero over $K$ (and hence satisfies the Hasse principle) as
soon as the number of variables $n$ is at least $9$.  This enables us
to handle an intersection $\cR$ of three quadrics over a totally
complex number field $K$, using the method of Leep \cite{leep}. If 
$n\geq 21$, then the projective variety
$Q_3=0$ automatically contains a linear space of dimension 8 defined
over $K$, since it must split off 9 hyperbolic planes. However the pair
of forms $Q_1,Q_2$ will have a common zero over $K$ in this linear
space, by the above mentioned 
result of Colliot-Th\'{e}l\`{e}ne, Sansuc and Swinnerton-Dyer, and
we deduce that any triple of quadratic forms in at least 21 variables, 
defined over a totally complex number field $K$, will have a common
zero in $K$.

However we wish to handle general number fields, and so we will assume that
the variety $\cR$ is nonsingular.
To be more precise, we shall require that the matrix
\[\left(\begin{array}{cc}\nabla Q_1(\x) \\ \nabla Q_2(\x) \\ \nabla Q_3(\x)
\end{array}\right)\]
has rank 3 for every point $[\x]$ in $\cR(\overline{K})$.  When this condition
holds we will say that $Q_1,Q_2,Q_3$ is a ``nonsingular system'' of quadratic
forms.  According to Lemma 3.2 of Browning and Heath-Brown
\cite{many}, this will ensure that $\cR$ is an absolutely
irreducible variety of codimension 3 and degree 8. For nonsingular
systems over $\Q$ one can apply the well-known result of Birch
\cite{birch}, which was proved using the Hardy--Littlewood circle
method.  Indeed Birch's work was generalized to arbitrary number
fields by Skinner \cite{skin}, who explicitly considers the question of
weak approximation. The outcome is that, for a nonsingular system,
$\cR$ will satisfy the Hasse principle and weak approximation provided
that $n>d+24$, where $d$ is the dimension of the ``Birch singular
locus''. In fact, for a nonsingular system of 3 quadratic forms the
Birch singular locus will have dimension at most 2, so that it
suffices to have $n\ge 27$.

Having described the relevant background we are ready to state our
principal result.
\begin{theorem}
Let $Q_1,Q_2,Q_3$ be a nonsingular system of quadratic forms in $n$
variables, defined over a number field $K$.  Then the variety
\[\cR:\; \cQ_1\cap\cQ_2\cap\cQ_3,\]
where $\cQ_i$ is the quadric $Q_i=0$, 
satisfies both the Hasse principle and weak approximation, as soon as
$n\ge 19$.
\end{theorem}
Thus we get both an improvement on the range $n\ge 21$ mentioned above
for the case of totally imaginary fields, and on the range $n\ge 27$
coming from the methods of Birch and Skinner. We should also observe
at this point that $\cR$ will have points over any completion $K_v$ 
at a finite place $v$, as soon as $n\ge 17$.  This is established by 
Heath-Brown \cite[page 138]{qbirs} when $K=\Q$, and the proof for
general number fields $K$ is completely analogous.  Thus as far as the
Hasse principle is concerned our theorem only requires solvability in
the real completions of $K$. In the same connection we mention that
$\cR$ has local points over $K_v$ for finite places over primes $p\ge
37$ as soon as $n\ge 13$, see Heath-Brown \cite[Corollary 1]{HBQs}.

Our basic strategy for proving the theorem will be to try and find
a linear space $L$ of dimension 7, defined over $K$ and lying in the
quadric hypersurface $\cQ_3$. If we can ensure that the variety
$\cQ_1\cap \cQ_2\cap L$ is nonsingular, we can then apply Theorem 1 of Heath-Brown
\cite{HB8}, as mentioned above. This approach will require us firstly
to establish the smoothness condition, and secondly to ensure that
$\cQ_1\cap \cQ_2\cap L$ has points everywhere locally. 

A comment is required about ``smoothness'' requirements, where we
shall adopt a significant abuse of terminology. When we say that the
variety $\cQ_1\cap \cQ_2\cap \cQ_3$ is nonsingular, for example, we mean that
the corresponding system of 3 quadratic forms is a nonsingular system,
as described above. It is of course quite possible for the variety 
$\cQ_1\cap \cQ_2\cap \cQ_3$ to be nonsingular in the usual sense, without
the corresponding system of forms being nonsingular. For example, if
all the forms $Q_i$ vanish identically then $\cQ_1\cap \cQ_2\cap
\cQ_3=\Proj^{n-1}$, which is nonsingular. When we say that $\cQ_1\cap
\cQ_2\cap L$ is nonsingular we will have in mind a system consisting of 2
quadratic forms and a (minimal) set of linear forms defining $L$.  A
variety can be defined in many ways via a system of forms, and the
reader will have to decide from the context what constitutes an
appropriate system to use.  Since all the relevant varieties will be
complete intersections there should be little difficulty with this.

\section{Geometric Considerations}

We begin by replacing the forms $Q_1,Q_2,Q_3$ by more convenient ones.
We will write $\bQ=(Q_1,Q_2,Q_3)$ for our triple of
quadratic forms, and proceed to consider linear combinations $\b{t}.\bQ
=t_1Q_1+t_2Q_2+t_3Q_3$.  The determinant $d_1(\b{t}):=\det(\b{t}.\bQ)$
is a form in $t_1,t_2,t_3$ of degree $n$. We also define the
determinant
\[\delta(X,Y;\b{t}^{(1)},\b{t}^{(2)}):=\det(X\b{t}^{(1)}.\b{Q}+Y\b{t}^{(2)}.\bQ)\]
and the discriminant
\[d_2(\b{t}^{(1)},\b{t}^{(2)}):=
{\rm Disc}(\delta(X,Y;\b{t}^{(1)},\b{t}^{(2)})).\]
Thus $\delta(X,Y;\b{t}^{(1)},\b{t}^{(2)})$ is a form of degree $n$ in
$X$ and $Y$, and
$d_2(\b{t}^{(1)},\b{t}^{(2)})$ is bihomogeneous in the entries of
$\b{t}^{(1)}$ and $\b{t}^{(2)}$.

It will be convenient to record some properties of nonsingular systems
of two quadratic forms.  The following lemma follows from Heath-Brown
and Pierce \cite[Proposition 2.1]{HBP}, for example. 
\begin{lemma}\label{hbp}
Let $R_1(x_1,\ldots,x_m)$ and $R_2(x_1,\ldots,x_m)$ be quadratic forms
over an algebraically closed field $k$ of characteristic zero, and
suppose that they constitute a nonsingular system, so that $\nabla R_1(\x)$
and $\nabla R_2(\x)$ are linearly independent for any non-zero $\x$
satisfying $R_1(\x)=R_2(\x)=0$.  Then every non-trivial linear
combination $aR_1+bR_2$ has rank at least $m-1$. 
Moreover $\det(XR_1+YR_2)$ is not
identically zero, and has
distinct linear factors over $k$. Conversely, this last condition is
equivalent to the statement that $R_1$ and $R_2$ form a nonsingular
system.  
\end{lemma}

We now use the following result of Aznar \cite[\S 2]{aznar}.
\begin{lemma}\label{lem:aznar}
Let  $V\subset \Proj^{n-1}$ be a nonsingular complete intersection  of
codimension $r$, which is defined over a field $k$ of characteristic
zero.  Then there is a system of generators $F_1,\ldots,F_r\in k[\x]$ of the
ideal of $V$, with 
\[\deg F_1\geq \cdots \geq \deg F_r,\]
such that the varieties
\[W_i: \quad F_1=\cdots=F_i=0, \quad (i\leq r),\]
are all non-singular. Moreover if $r<n-1$ then $V$ will be irreducible 
with degree equal to the product of the degrees of the $F_i$.
\end{lemma}
The last statement in the lemma follows from Heath-Brown and Browning
\cite[Lemma 3.2]{many}.
When $V$ is originally defined by $r$ forms $G_j$ of the same
degree $d$, the proof of Aznar's result
shows that one can take the $F_i$ to be linear
combinations of the $G_i$. 

In our case Lemma \ref{lem:aznar} implies in particular 
that there are linear combinations
$\b{m}^{(1)}.\bQ$ and $\b{m}^{(2)}.\bQ$, with
$\b{m}^{(1)},\b{m}^{(2)}\in \overline{K}^3$, such that both the hypersurface
$\b{m}^{(1)}.\bQ=0$ and the intersection
$\b{m}^{(1)}.\bQ=\b{m}^{(2)}.\bQ=0$ are nonsingular.  In particular we
will have $d_1(\b{m}^{(1)})\not=0$, so that the form $d_1(\b{t})$ does
not vanish identically. Moreover it follows from Lemma \ref{hbp} 
that the intersection
$\b{m}^{(1)}.\bQ=\b{m}^{(2)}.\bQ=0$ is nonsingular if and only if 
$d_2(\b{m}^{(1)},\b{m}^{(2)})$ is non-zero, and so we deduce that the form
$d_2(\b{t}^{(1)},\b{t}^{(2)})$ does not vanish identically. 

We therefore see that if we choose any three vectors $\b{m}^{(1)},
\b{m}^{(2)}, \b{m}^{(3)}\in K^3$ such that none of
\[d_1(\b{m}^{(3)}),\;\; d_2(\b{m}^{(1)},\b{m}^{(2)}),\;\;
d_2(\b{m}^{(1)},\b{m}^{(3)}),\;\;\mbox{or}\;\;
\det\big(\b{m}^{(1)}|\b{m}^{(2)}|\b{m}^{(3)}\big),\]
vanishes,
then the three forms $Q_i'(\x)=\b{m}^{(i)}.\bQ(\x)$, for $i=1,2,3$ are
defined over $K$ and 
generate the same linear system as do $Q_1,Q_2,Q_3$.  Moreover $Q_3'=0$
is nonsingular, as are the varieties $Q_1'=Q_2'=0$ and $Q_1'=Q_3'=0$. 
Thus without loss of generality we will assume that our
original forms satisfy these conditions.  

We will also require $\cQ_3$ to contain suitable linear spaces defined
over the real completions of $K$.  Suppose as above
that $Q_1$ and $Q_3$ form a
nonsingular system, and that $K_v$ is real.  Then the argument of 
Heath-Brown \cite[Lemma 12.1]{HB8} shows that there is a real $\theta_v$ for
which $(\cos\theta_v)Q_1+(\sin\theta_v)Q_3=0$ is nonsingular and contains
a linear space of dimension at least $(n-4)/2$ over $K_v$.  (The
argument of \cite{HB8} does not explicitly produce a nonsingular
quadratic form, but since the functions $n_+$ and $n_-$ are everywhere
locally minimal one can change $\theta_v$ slightly, if necessary.)
By weak
approximation in $\Proj^1(K)$ we deduce that there exists $c\in K$ for
which $cQ_1+Q_3=0$ is also nonsingular and contains
linear spaces of dimension at least $(n-4)/2$ over each real
completion $K_v$.  We now replace $Q_3$ by $cQ_1+Q_3$ so that $\cQ_3$ is
nonsingular and has linear spaces over
each real completion, of dimension at least $(n-4)/2$.  
For finite places $v$ it is automatic that $\cQ_3$
contains linear spaces over $K_v$, of dimension at least $(n-5)/2$.  We can
therefore conclude that $\cQ_3$ has a linear space over $K$, with
dimension at least $(n-5)/2$. The existence of a single such linear space
is enough to ensure that there is one through every $K$-point of
$\cQ_3$. An alternative way to express the above facts is to say that
$Q_3$ splits off at least $(n-5)/2$ hyperbolic planes over $K$.

From now on we will fix the forms $Q_1$, $Q_2$ and $Q_3$.  We remind
the reader that the varieties $\cR$, $\cQ_1\cap\cQ_2$ and $\cQ_3$ are
all nonsingular, and that we have arranged that $\cQ_3$ contains a linear
space, defined over $K$, of dimension at least $(n-5)/2$.

For the remainder of this section we work over an arbitrary 
algebraically closed field $k$ of characteristic zero.
After the above preliminary manoeuvres we are ready to prove our first
key result.

It will be convenient to define
\[F_t=F_t(\cQ_3):=\{L\in\mathbb{G}(t,n-1):\, L\subset \cQ_3\}.\]
\begin{lemma}\label{existL}
Suppose that the forms $Q_1,Q_2,Q_3$ are such that the varieties
$\cR=\cQ_1\cap\cQ_2\cap\cQ_3$, $\cQ_1\cap\cQ_2$, and $\cQ_3$ are all
nonsingular.
Then for every integer $t$ in the range $3\le 
t\le(n-5)/2$ there is an $L\in F_t$,
such that $L\cap\cQ_1\cap\cQ_2$ is nonsingular.
\end{lemma}
 
We should perhaps be more specific as to this last condition.  If $L$
has codimension $c$ say, and is given by linear equations
$\ell_1(\x)=\ldots=\ell_c(\x)=0$, then according to our convention,
$L\cap\cQ_1\cap\cQ_2$ is nonsingular if the partial derivatives of
$Q_1,Q_2$ and $\ell_1,\ldots,\ell_c$ are linearly independent at any
non-zero vector $\x\in k^n$ such that 
\[Q_1(\x)=Q_2(\x)=\ell_1(\x)=\ldots=\ell_c(\x)=0.\]
Thus for example, if $L\subset \cQ_3$ is a line disjoint from 
$\cQ_1\cap\cQ_2$, then $L$ fulfils the condition vacuously.

For the proof we will need some information about $t$-planes lying in
quadric hypersurfaces. We begin by introducing some notation.
Let $Q(x_1,\ldots,x_n)$ be a quadratic form of rank $r$, over 
$k$, and let $\cQ$
be the quadric $Q=0$, with dimension $n-2$.  Let $F(n,r,t)$ be the 
Fano variety of $t$-planes in $\cQ$, and let $F(n,r,t;P)$ be the
subvariety of such planes passing through a given point $P\in \cQ$.
The variety $F(n,r,t)$ will be non-empty when $t\le n-r/2-1$. 
Write $d_0(n,t,r)=\dim F(n,t,r)$, which will be independent of the
particular quadratic form $Q$.  Similarly let
$\dim F(n,t,r;P)=d_1(n,t,r)$ for $P$ a smooth point of $\cQ$.
These dimensions are given by the following lemma.

\begin{lemma}\label{fano}
We have the following statements.
\begin{enumerate}
\item[(i)] $d_1(n,t,r)=d_0(n-2,t-1,r-2)$ if $t\ge 1$ and $r\ge 2$.
\item[(ii)] $d_0(n,t,r)=n-2-t+d_0(n-2,t-1,r-2)$ if $t\ge 1$ and 
$2\le r\le 2n-2t-2$.
\item[(iii)] $d_0(n,t,r)=(t+1)(n-2-3t/2)$ if $2t+2\le r\le 2n-2t-2$.
\end{enumerate}
\end{lemma}

We take $P=[\b{e}_1]$ and extend to a basis
$\b{e}_1,\ldots,\b{e}_n$ of $K^n$. Then $Q$ takes the shape
$x_1L(x_2,\ldots,x_n)+Q'(x_2,\ldots,x_n)$ with respect to this basis,
with $L\not=0$. A further change of basis simplifies this to $x_1x_2+
Q''(x_3,\ldots,x_n)$. Here $Q''$ will have rank $r-2$.  One then sees
that $t$-planes in $\cQ$ containing $\b{e}_1$ correspond
to $(t-1)$-planes in $Q''=0$, and the result (i) follows. 

For part (ii) we consider the incidence correspondence
\[I=\{(P,L):P\in L\in F(n,t,r)\}.\]
The projection $\pi_2$ onto the second factor takes $I$ onto
$F(n,t,r)$, and each fibre has dimension $t$, so that
$d_0(n,t,r)=\dim(I)-t$. The condition $r\le 2n-2t-2$ ensures that
$d_0(n-2,t-1,r-2)\ge 0$, whence $d_1(n,t,r)\ge 0$ by part (i).  Thus,
for the projection $\pi_1$ onto the first factor, we see that
$\pi_1(I)$ contains every smooth point of $\cQ$.  Hence
$\pi_1(I)=\cQ$, and $\pi_1^{-1}(P)$ will have dimension $d_1(n,t,r)$ for
smooth points $P$.  Thus $\dim(I)=\dim(\cQ)+d_1(n,t,r)$, and 
the result follows from part (i).

Finally, part (iii) follows from part (ii) by induction on $t$, the
result being clearly true for $t=0$.

We can now move to the proof of Lemma \ref{existL}
By Lemma \ref{fano} we have
\[\dim(F_t)=d_0(n,t,n)=(t+1)(n-2-3t/2)\]
if $2t+2\le n$.
We proceed to consider
the variety $F^{\dagger}$ defined to be the set of $t$-planes $L\in F$ for
which $L\cap\cQ_1\cap\cQ_2$ is singular.  Let
\[I:=\{(L,[\x],[\b{t}])\in F_t\times\cR\times\Proj^1:\,
[\x]\in L,\,\b{y}^T(t_1Q_1+t_2Q_2)\x=0\,
\forall [\b{y}]\in L\}.\]
If $L\cap\cQ_1\cap\cQ_2$ has a singular singular point $\x$, then 
$\x\in L\cap\cR$ and there is some $[\b{t}]\in\Proj^1$ such that
$\b{y}^T(t_1Q_1+t_2Q_2)\x=0$ for every $[\b{y}]\in L$.
Thus if $\pi_1$ is the projection from $I$ onto its first factor,
one has $\pi_1(I)=F^{\dagger}$. It follows that $\dim(F^{\dagger})\le\dim(I)$.

Now consider the projection $\pi_{2,3}$ onto the second and third
factors. The fibre above the pair $([\x],[\b{t}])$ will be
\beql{fib}
\{L\in F_t:\,[\x]\in L,\,\b{y}^T(t_1Q_1+t_2Q_2)\x=0\,\forall [\b{y}]\in L\}.
\eeq
We write $Q=t_1Q_1+t_2Q_2$ and
\[H(\x,\b{t})=\{[\b{y}]:\,\b{y}^TQ\x=0\}\]
for convenience. Thus the fibre (\ref{fib}) may be written as
\[\{L\in F_t:\, L\subseteq H(\x,\b{t}),\,[\x]\in L\}.\]
Since $[\x]\in\cR$, and $\cR$ is nonsingular, we must have 
$Q\x\not=\b{0}$, so that $H(\x,\b{t})$ is a hyperplane. Hence
$H(\x,\b{t})$ intersects $\cQ_3$
to produce a quadric hypersurface $\cQ'$ in
$\Proj^{n-2}$, with rank, $r$ say, at least $n-2$. The point $[\x]$
must be a smooth point of $\cQ'$, since otherwise $Q\x$ and $Q_3\x$
would be proportional, contradicting the fact that $[\x]$ is a smooth
point of $\cR$.

It follows that
the dimension of the fibre (\ref{fib}) will be $d_1(n-1,t,r)$.
According to Lemma \ref{fano} we have
\[d_1(n-1,t,r)= d_0(n-3,t-1,r-2)=t\big(n-5-3(t-1)/2\big)\]
if $2(t-1)+2\le r-2\le 2(n-3)-2(t-1)-2$.  
Since we are assuming that $t\le (n-5)/2$ the
required condition on $r$ certainly holds. We therefore deduce that
\begin{eqnarray*}
\dim(F^{\dagger})\leq\dim(I)&=&\dim(\cR\times\Proj^1)
+\dim\pi_{2,3}^{-1}([\x],[\b{t}])\\
&=& \big((n-4)+1\big)+t\big(n-5-3(t-1)/2\big)\\
&=& (t+1)\big(n-2-3t/2)-1\\
&=& \dim(F_t)-1
\end{eqnarray*}
so that $F^{\dagger}$ must be a proper subvariety of $F_t$.
This completes the proof of Lemma \ref{existL}.
\bigskip

We turn now to our second key result.
\begin{lemma}\label{existL2}
Suppose that the forms $Q_1,Q_2,Q_3$ are such that the varieties
$\cR=\cQ_1\cap\cQ_2\cap\cQ_3$, $\cQ_1\cap\cQ_2$, and $\cQ_3$ are all
nonsingular.
Then for every non-negative integer $t\le(n-5)/2$ there is an
$L\in F_t$ such that $(Q_3L)^{\perp}\cap\cR$ is nonsingular.
\end{lemma}
Here we define
\[(Q_3L)^{\perp}:=\{[\b{y}]:\,\b{y}^TQ_3\x=0,\,\forall [\x]\in L\}.\]

We prove this by induction on $t$.
For the base case $t=0$ the space $L$ is a single point $P$, say. It then
suffices that $P\in \cQ_3$ and $P\not\in
Q_3^{-1}(\cR^*)$, where $\cR^*$ is the dual variety to $\cR$, and
$Q_3^{-1}(\cR^*)$ is its image under the linear map $Q_3^{-1}$.
There will always
be a suitable point $P$ if $\cQ_3\not\subseteq Q_3^{-1}(\cR^*)$. However $\cR^*$
is a proper subvariety of $\Proj^{n-1}$.  We claim that it cannot be a
nonsingular quadric, whence we cannot have $\cQ_3=Q_3^{-1}(\cR^*)$. To prove the
claim we merely observe that if $\cR^*=\cQ$, say, then
$\cR=\cR^{**}=\cQ^*$. However $\cQ^*$ is itself a nonsingular quadric,
whereas $\cR$ has degree 8, by Lemma \ref{lem:aznar}. The lemma then
follows in the case of dimension zero.

To establish the induction step we suppose we have a suitable space
$L$ of dimension $t$, and look for an appropriate $L'$ of dimension
$t+1$. Indeed we shall restrict our attention to linear spaces satisfying
\[L\subset L'\subseteq\cQ_3\cap (Q_3L)^{\perp}.\]
This requirement on $L'$ allows us to restrict all our varieties to the
subspace $(Q_3L)^{\perp}$. We write $\cQ_i'=\cQ_i\cap(Q_3L)^{\perp}$
for $i=1,2,3$, and $\cR'=\cR\cap(Q_3L)^{\perp}$. Thus $\cR'$ is nonsingular.
More concretely, we choose
$[\b{e}_0],\ldots,[\b{e}_t]$ spanning $L$, and extend these to a set of points
$[\b{e}_0],\ldots,[\b{e}_{n-t-2}]$ spanning $(Q_3L)^{\perp}$, and then to a set
$[\b{e}_0],\ldots,[\b{e}_{n-1}]$ spanning $\Proj^{n-1}$. We proceed to
write
\[Q_i'(x_1,\ldots,x_{n-t-1})=Q_i(\sum_{j=1}^{n-t-1}x_j\b{e}_{j-1})
\;\;\;(i=1,2,3)\]
so that these quadratic forms correspond to the quadrics $\cQ_i'$,
seen as varieties in $\Proj^{n-t-2}$. We note in particular that
$Q_3'$ is singular, with rank $n-2t-2$. Indeed, as a subvariety of 
$\Proj^{n-t-2}$ the quadric $\cQ_3'$ hypersurface will be a cone with
vertex set $L\subset\Proj^{n-t-2}$.

We construct
$L'$ as $<L,[\x]>$ with $[\x]\in\cQ_3'$ but $[\x]\not\in L$.  
This will ensure that $L'\subseteq\cQ_3'$ and that $\dim(L')=t+1$. 
Moreover we will have
$(Q_3L')^{\perp}=(Q_3L)^{\perp}\cap(Q_3\x)^{\perp}$, so that
$(Q_3L')^{\perp}\cap\cR'$ will be nonsingular provided that $[Q_3\x]$
is not in $(\cR')^*$. Since $\dim(L)=t<\dim(\cQ_3')=n-t-3$, the generic
$[\x]\in\cQ'_3$ will satisfy $[\x]\not\in L$.  Moreover it will also
satisfy $[Q_3\x]\not \in (\cR')^*$ by a similar reasoning to that given
in the dimension zero case above.  Specifically, $(\cR')^*$ is a
proper subvariety of $\Proj^{n-t-2}$, so the only situation to rule
out is that in which it is equal to $\cQ_3'$. However if $Q_3'$ has
rank $n-2t-2$ then the dual $(\cQ_3')^*$ will be a quadric of dimension 
$n-2t-4$, which cannot possibly be the variety $\cR'$, since the
latter will have degree 8.

It follows that we can take $L'=<L,[\x]>$ for a generic
$[\x]\in\cQ_3'$. This completes the induction step and so establishes
the lemma.
\bigskip

Following on from Lemmas \ref{existL} and \ref{existL2} we make
the following definitions.
\begin{definition}
Let $Q_1,Q_2,Q_3$ be a nonsingular system of quadratic forms in $n$
variables as above, and let
$L\subseteq\Proj^{n-1}$ be a linear space contained in
the quadric $\cQ_3$.  Then we say that $L$ is ``admissible'' if and
only if both $L\cap\cQ_1\cap\cQ_2$ and $(Q_3L)^{\perp}\cap\cR$ are nonsingular.
If $L$ has dimension $t$ in the range $3\le t\le 7$ 
we say that $L$ is ``chain-admissible'' if
there exist admissible linear spaces $L=L_t\subset
L_{t+1}\subset\ldots\subset L_7$ with $\dim(L_i)=i$ for $t\le i\le 7$.
\end{definition}

We then have the following result.
\begin{lemma}\label{existL3}
For each $t\in[3,7]$ there is a Zariski-closed proper subset $Z_t\subset
F_t$, such that $L\in F_t$ is chain-admissible if and only if
$L\not\in Z_t$.
\end{lemma}

It is clear that the $L\in F_t$ for which $L\cap\cQ_1\cap\cQ_2$ is
singular form a closed subset $A_t$ say, of $F_t$.  Similarly those
for which $(Q_3L)^{\perp}\cap\cR$ is singular form a closed subset
$B_t$ say. By Lemmas \ref{existL} and \ref{existL2} these are
proper closed subsets of $F_t$.  Moreover, since $F_t$ is absolutely
irreducible, the union $A_t\cup B_t=C_t$, say, is also a closed proper
subset of $F_t$.  A linear subspace $L\in F_t$ is admissible if and
only if $L\not\in C_t$. We proceed to prove Lemma \ref{existL3}
by downwards induction, and it is clear that we can take $Z_7=C_7$.

For $t<7$ if $L\in F_t$ fails to be chain-admissible then either $L\in
C_t$ (because $L$ is not itself admissible) or $L'$ fails to be
chain-admissible for every $L'\in F_{t+1}$ containing $L$.  This
latter case holds precisely when $L'\in Z_{t+1}$ for every such $L'$. Write
$D_t$ for the set of $L\in F_t$ with the property that $L'\in Z_{t+1}$
for every $L'\in F_{t+1}$ containing $L$.  We claim that $D_t$ is a
proper closed subset of $F_t$.  Once this is established the induction
step of the proof is completed by taking $Z_t=C_t\cup D_t$.

To handle $D_t$ we consider 
\[I:=\{(L,L')\in F_t\times Z_{t+1}:\, L\subset L'\}.\]
For the projection onto the first factor we have
\[\dim(\pi_1^{-1}(L))=\dim\{L'\in Z_{t+1}:\,L\subset L'\}\leq n-2t-4,\]
with equality exactly when $L'\in Z_{t+1}$ for every linear space $L'\in
F_{t+1}$ containing $L$. Thus $D_t$ is the set of $L$ for which
$\dim(\pi_1^{-1}(L))$ is maximal, whence $D_t$ is
Zariski-closed. Moreover we have
\[\dim(D_t)\le \dim(I)-(n-2t-4).\]
On the other hand, for the projection onto the second factor we have
\[\dim(\pi_2^{-1}(L'))= t+1,\]
whence
\[\dim(I)=\dim(Z_{t+1})+t+1.\]
Since $Z_{t+1}$ is a proper subset of $F_{t+1}$ by the downward
induction hypothesis we deduce that
\begin{eqnarray*}
\dim(D_t)&\le& \dim(I)-(n-2t-4)\\
&=& \dim(Z_{t+1})+t+1-(n-2t-4)\\
&<& \dim(F_{t+1})-n+3t+5\\
&=& (t+2)(n-2-3(t+1)/2)-n+3t+5\\
&=& (t+1)(n-2-3t/2) \\
&=& \dim(F_t).
\end{eqnarray*}
Here we have used Lemma \ref{fano} to compute
\[\dim(F_{t+1})=d_0(n,t+1,n)=(t+2)(n-2-3(t+1)/2)\]
and
\[\dim(F_t)=d_0(n,t,n)=(t+1)(n-2-3t/2).\]
The above calculation shows that $\dim(D_t)<\dim(F_t)$, so that $D_t$
is a proper subset of $F_t$, as required.

We conclude this section with two easy results in a similar vein.
\begin{lemma}\label{R3chain}
There is a Zariski-closed proper subset $\cR_0\subset\cR$
such that every $P\in\cR-\cR_0$ is contained in a chain-admissible
linear space $L\in F_3-Z_3$.
\end{lemma}

Let $\cR_0$ be the set of points $P\in\cR$ such that $L\in Z_3$ for every
$L\in F_3$ which contains $P$.  We need to show that $\cR_0$ is a
proper closed subset of $\cR$.  It follows from Lemma \ref{existL3}
that there is at least one $L\in F_3-Z_3$.  Since
$\dim(L)+\dim(\cR)=3+(n-4)=n-1$ it follows that $L\cap\cR$ is
non-empty, containing $P$ say.  Then $P\in \cR-\cR_0$, so that $\cR_0$
is a proper subset of $\cR$.

To show that $\cR_0$ is Zariski-closed we consider
\[I=\{(P,L)\in \cR\times Z_3:\, P\in L\}.\]
For the projection onto the first factor we have
\[\dim(\pi_1^{-1}(P))\leq d_1(n,3,n),\]
in the notation of Lemma \ref{fano}, with equality exactly when $L\in
Z_3$ for every $L\in F_3$ which contains $P$.  Thus $\cR_0$ is the set
of $P$ for which $\dim(\pi_1^{-1}(P))$ is maximal, and it follows that
$\cR_0$ is Zariski-closed, as claimed.

\begin{lemma}\label{good4}
Given $L\in F_3-Z_3$, there is a Zariski-closed proper subset $\cR_1$
of $\cR\cap\M$, such that $<L,P>\in F_4-Z_4$ for every
$P\in\cR\cap\M-\cR_1$. Similarly, given $L\in F_t-Z_t$, for some
$t\in[3,6]$, there is a Zariski-closed proper subset $E_t$ of
$\cQ_3\cap\M$, such that $<L,P>\in F_{t+1}-Z_{t+1}$ for every
$P\in\cQ_3\cap\M-E_t$. 
\end{lemma}

To prove the first part of the lemma we let
\[I:=\{(P,L')\in(\cR\cap\M)\times Z_4:\, P\in L',\,L\subset L'\},\]
and take $\cR_1=\pi_1(I)$, so that $\cR_1\subseteq\cR\cap\M$ 
is clearly Zariski-closed.  Since $L\in F_3-Z_3$ is chain-admissible
there is at least one linear space $L_0\in F_4-Z_4$ containing
$L$. Then $L_0\cap\cR=L_0\cap\cQ_1\cap\cQ_2$ has dimension $2$, since
$L_0$ is admissible, and similarly $L\cap\cR=L\cap\cQ_1\cap\cQ_2$ 
has dimension $1$. We may therefore find a point $P\in
L_0\cap\cR-L\cap\cR$. We claim that $P\in\cR\cap\M-\cR_1$, which 
shows that $\cR_1$ is a proper subset of $\cR\cap\M$, and so proves the
first part of the lemma.  Since $L\subset L_0\subset\cQ_3$ it follows that
$L_0\subseteq\M$, whence $P\in L_0\cap\cR\subseteq \M\cap\cR$.  On the
other hand if we had $P\in\cR_1$ there would be a linear space $L'$
such that $(P,L')\in I$.  Then we would have $P\in L'$ and $L\subset
L'$, and therefore $L'=<L,P>$, since $P\not\in L$ by our choice of
$P$. However the same reasoning shows that $L_0$ is also equal to
$<L,P>$, so that $L'=L_0$.  This gives us a contradiction, since
$L'\in Z_4$ while $L_0\in F_4-Z_4$.

Turning to the second part of the lemma, we consider
\[I:=\{(P,L')\in(\cQ_3\cap\M)\times Z_{t+1}:\, P\in L',\,L\subset L'\},\]
and take $E_t=\pi_1(I)$, so that $E_t\subseteq\cQ_3\cap\M$ 
is clearly Zariski-closed.  Since $L\in F_t-Z_t$ is chain-admissible
there is at least one linear space $L_0\in F_{t+1}-Z_{t+1}$ containing
$L$. We claim that $P\in\cQ_3\cap\M-E_t$ for any point $P\in L_0-L$,
which will show that $E_t$ is a proper subset of $\cQ_3\cap\M$.
Since $L\subset L_0\subset\cQ_3$ it follows that
$L_0\subseteq\M$, whence $P\in\cQ_3\cap\M$.  On the
other hand if we had $P\in E_t$ there would be a linear space $L'$
such that $(P,L')\in I$.  Then we would have $P\in L'$ and $L\subset
L'$, and therefore $L'=<L,P>$, since $P\not\in L$ by our choice of
$P$. However the same reasoning shows that $L_0$ is also equal to
$<L,P>$, so that $L'=L_0$.  This gives us a contradiction, since
$L'\in Z_{t+1}$ while $L_0\in F_{t+1}-Z_{t+1}$.

\section{Global 3-planes in $\cQ_3$}

In this section we shall make repeated use of three key principles.
The first of these is the fact that we have weak approximation on
quadrics.  The second is that if $V$ is an absolutely irreducible
projective variety defined over the number field $K$, with a smooth 
point $P_v$ in some completion $K_v$, then the $K_v$-points of $V$ are
Zariski-dense in any given neighbourhood of $P_v$. (See Browning,
Dietmann and Heath-Brown
\cite[Lemma 3.4]{BDHB}, for example. The proof is an application of
the implicit function theorem.) 

The third general principle is embodied in the following lemma.
\begin{lemma}\label{approxV}
Let $V$ be a projective 
algebraic variety defined over a completion $K_v$ of $K$
by equations
\[f_1(x_0,\ldots,x_m)=\ldots=f_r(x_0,\ldots,x_m)=0,\]
and suppose that $P$ is $K_v$-point on $V$ at which the vectors
$\nabla f_1,\ldots,\nabla f_r$ are linearly independent. Suppose we
are given varieties $V^{(j)}$ defined over $K_v$ by equations
\[f_1^{(j)}(x_0,\ldots,x_m)=\ldots=f_r^{(j)}(x_0,\ldots,x_m)=0,\]
of bounded degree,
in which $f_i^{(j)}\rightarrow f_i$ (under the metric induced from
$K_v$) as $j\rightarrow\infty$.  Then for sufficiently large $j$ there
are $K_v$-points $P^{(j)}$ on $V^{(j)}$ such that $P^{(j)}\rightarrow
P$ as $j\rightarrow\infty$. Moreover $P^{(j)}$ will be a nonsingular
point of $V^{(j)}$, in the sense above.
\end{lemma}

For the proof we suppose firstly that $v$ is a finite place. We let
$P=[\b{t}]$ say, and we rescale $\b{t}$ and the polynomials $f_j$ so
as to have $v$-adic integer coefficients. Since $f_i^{(j)}\rightarrow f_i$ it
follows that $f_i^{(j)}$ also has integral coefficients if $j$ is large
enough. By hypothesis, the matrix
formed from the rows $\nabla f_1(\b{t}),\ldots,\nabla f_r(\b{t})$ has
rank $r$. We may therefore suppose without loss of generality that the
determinant, $\Delta$ say, of the first $r$ columns is non-zero. Let
$\Delta_j$ be the corresponding determinant formed from
$\nabla f_1^{(j)}(\b{t}),\ldots,\nabla f_r^{(j)}(\b{t})$. Thus
$\Delta_j\rightarrow\Delta$, so that $|\Delta_j|_v=|\Delta|_v\not=0$ if $j$ is
large enough.  

We also set
\[\delta_j=\max\{|f_1^{(j)}(\b{t})|_v,\ldots,|f_r^{(j)}(\b{t})|_v\},\]
and note that $\delta_j$ tends to
\[\max\{|f_1(\b{t})|_v,\ldots,|f_r(\b{t})|_v\}=0,\]
since $P\in V$. Thus if $j$ is large enough we will have
$\delta_j<|\Delta|_v^2$. This condition allows us to use Hensel's Lemma,
which provides a point $P^{(j)}=[\b{t}^{(j)}]$ on $V^{(j)}$, with 
\beql{henb}
\max_{i}|t_i^{(j)}-t_i|_v\le \delta_j/|\Delta|_v.
\eeq
It follows that $P^{(j)}$ tends to $P$ as required. Since
$\Delta_j\not=0$ for large $j$ the nonsingularity condition also holds.

When $v$ is an infinite place we use a completely analogous argument,
replacing Hensel's Lemma by Newton Approximation. We begin by
normalizing so that the entries of $\b{t}$, and the coefficients of
the $f_i$, all have modulus at most 1.  The condition 
$\delta_j<|\Delta|^2$ has to be replaced by $\delta_j<C|\Delta|^2$
with a constant $C$ depending on $m,r$ and the degrees of the
polynomials involved. Similarly, the bound in (\ref{henb}) becomes
$C'\delta_j/|\Delta|$ with a corresponding constant $C'$.  With these
changes the proof goes through as before.
\bigskip

For our theorem we will assume that we are given local points
$[\x_v]\in \cR(K_v)$ for every place $v$ of $K$.  We will also be
given a finite set of places $S$ and a (small) positive $\ep$, and our
challenge will be to find a point $[\x]\in\cR(K)$ such that
\beql{appr}
\left|\x-\x_v\right|_v<\ep\;\;\;\mbox{for all}\; v\in S.
\eeq
Without loss of generality we will include all infinite places in $S$,
as well as all finite places above rational primes up to 37.
In particular, $S$ will be non-empty.
From now on we will assume that the number $n$ of variables in our
quadratic forms satisfies $n\ge 19$. 

The variety $\cR$ is nonsingular, and for each
$v\in S$ the $K_v$-points of $\cR$ are therefore Zariski-dense in
every neighbourhood of $[\x_v]$, by the second principle above.  
It follows that if $\cR_0$ is as in Lemma \ref{R3chain} then
there is a point
$[\x_v']\in \cR-\cR_0$, defined over $K_v$, in the neighbourhood
$\left|\x_v-\x_v'\right|_v<\ep/2$. Thus it suffices to find a
$K$-point of $\cR$ with
\[\left|\x-\x_v'\right|_v<\ep/2\;\;\;\mbox{for all}\; v\in S,\]
where now there is a chain-admissible 3-plane $L_v\in F_3$ through 
$[\x_v']$.  We therefore change our
notation, replacing $\x_v'$ by $\x_v$ and $\ep$ by $\ep/2$ so that is
still suffices to work with the condition (\ref{appr}). 
Note that $L_v$ may be
defined over $\overline{K_v}$ rather than $K_v$. None the less
the existence of a single chain-admissible $L_v$ shows that the generic
3-plane $L\subset\cQ_3$ through $[\x_v]$ is also chain-admissible.
Thus we can in fact assume that $L_v$ is defined over
$K_v$.

Our plan now is to produce a sequence of 3-planes $L^{(m)}$ defined
over $K$, which approximate $L_v$ for each $v\in S$, in the following
sense. 
\begin{lemma}\label{Lappr}
For each $v\in S$ let $L_v\in F_3-Z_3$ be a chain-admissible
3-plane defined over $K_v$, and suppose that
$[\b{e}_{0,v}],\ldots,[\b{e}_{3,v}]$ is a basis of $L_v$.
Then there are sequences of $K$-points
$[\b{e}_{0}^{(m)}],\ldots,[\b{e}_{3}^{(m)}]$ spanning chain-admissible
3-planes 
$L^{(m)}\subset\cQ_3$, such that
\beql{lim!}
\lim_{m\rightarrow\infty}\b{e}_i^{(m)}=\b{e}_{i,v}\;\;\;\mbox{for
}0\le i\le 3,
\eeq
for every $v\in S$.
\end{lemma}

(We teach our students that a sequence can have at most one limit!
There is of course some abuse of notation above, and strictly speaking
we should have said that $\iota_v(\b{e}_i^{(m)})$ tends to
$\b{e}_{i,v}$, where $\iota_v$ is the embedding of $K^n$ into $K_v^n$.)
In effect the result is saying that the Fano variety of $3$-planes in
$\cQ_3$ satisfies weak approximation.
We will produce the sequences $[\b{e}_t^{(m)}]$ by induction on $t$, 
the case $t=0$ merely being an instance of weak approximation on $\cQ_3$. 

We therefore consider the induction step, and suppose we already have
suitable sequences of vectors
$\b{e}_{0}^{(m)},\ldots,\b{e}_{t-1}^{(m)}$. Thus the conditions
required for $\b{e}_t^{(m)}$ are that
\[[\b{e}_t^{(m)}]\in V^{(m)}:=\{[\x]\in\cQ_3: \x^TQ_3\b{e}_i^{(m)}=0
\mbox{ for }0\le i\le t-1\},\]
and that
\beql{lim}
\lim_{m\rightarrow\infty}\b{e}_t^{(m)}=\b{e}_{t,v}
\eeq
for every $v\in S$.  Notice that (\ref{lim!}) (for $i\le t-1$)
and (\ref{lim}) automatically ensure
that the vectors $\b{e}_{0}^{(m)},\ldots,\b{e}_t^{(m)}$ are linearly
independent, if $m$ is large enough.
The variety $V^{(m)}$ is a quadric of rank at least
$n-2t\ge 5$, so that it must have smooth points over $K_v$ for every finite
place $v$. For every $v\in S$, and in particular for every infinite
place, the varieties $V^{(m)}$ are approximations to
\[V_v:=\{[\x]\in\cQ_3: \x^TQ_3\b{e}_{i,v}=0
\mbox{ for }0\le i\le t-1\},\]
in the sense given by Lemma \ref{approxV}.  Moreover $[\b{e}_{t,v}]$
is a nonsingular point on $V_v$, in the sense of the lemma. It
follows that for large enough $m$ the variety $V^{(m)}$ has a smooth point
$[\b{f}_v^{(m)}]$ for every $v\in S$, such that 
\[\lim_{m\rightarrow\infty}\b{f}_v^{(m)}=\b{e}_{t,v}.\]
In particular $V^{(m)}$ has points everywhere locally, and so has a
$K$-point, by the Hasse principle.

We can now complete the induction step.  Given $\eta>0$ we choose
$m_0(\eta)$ so that
\[|\b{f}_v^{(m)}-\b{e}_{t,v}|_v<\eta/2\]
for all $v\in S$, and all $m\ge m_0(\eta)$. Moreover we can use weak
approximation on $V^{(m)}$ to find points $\b{e}_t^{(m)}$ on $V^{(m)}(K)$ such that
\[|\b{e}_t^{(m)}-\b{f}_v^{(m)}|_v<\eta/2\]
for all $v\in S$, and all $m\ge m_0(\eta)$. Then for all $v\in S$ and
all $m\ge m_0(\eta)$ we will have
\[|\b{e}_t^{(m)}-\b{e}_{t,v}|_v<\eta,\]
whence (\ref{lim}) holds, as required.
Finally, since $L_v\not\in Z_3$ we will have $L^{(m)}\not\in Z_3$ for
large enough $m$, by
continuity, so that $L^{(m)}$ is also chain-admissible.

\section{Completion of the Argument}

Our strategy now is to consider the intersection of $\cQ_1\cap\cQ_2$
with $L^{(m)}$. We begin with the following result.
\begin{lemma}\label{3p}
Let $\ep>0$ be given.  If $m$ is large enough, for every 
$v\in S$ there will be a $K_v$-point $[\b{y}_v^{(m)}]
\in\cQ_1\cap\cQ_2\cap L^{(m)}$ such that
\[|\b{y}_v^{(m)}-\x_v|_v<\ep/2.\]
\end{lemma}

This is a further application of Lemma \ref{approxV}.  The varieties
\[V^{(m)}:=\cQ_1\cap\cQ_2\cap L^{(m)}\]
are approximations to 
$V:=\cQ_1\cap\cQ_2\cap L_v$, in the sense of the lemma, and $P=[\x_v]$
lies on $V$, since we have both $P\in\cR$ and $P\in L_v$. Moreover 
$L_v$ is chain-admissible, and hence in particular is admissible, whence
$V$ is nonsingular. The lemma therefore produces appropriate points
$[\b{y}_v^{(m)}]$. 
\bigskip

We now fix a suitable $m$ in Lemma \ref{3p}, and write 
$L=L^{(m)}$ and $\b{y}_v=\b{y}_v^{(m)}$ accordingly.  Thus $L$ is a
chain-admissible 3-plane, defined over $K$. Moreover for each
$v\in S$ we have $[\b{y}_v]\in\cQ_1\cap\cQ_2\cap L$.  Finally, to
prove our theorem it will suffice to find a point $[\x]\in\cR(K)$ with
\[|\b{y}_v-\x|_v<\ep/2\]
for each $v\in S$.

For each $v\in S$ the variety $\cQ_1\cap\cQ_2\cap L$ has a local point, namely
$[\b{y}_v]$.  Moreover it is nonsingular, since $L$ is chain-admissible.
Thus there are local points at all but
finitely many places.  Let $T$ be the set of places for which there
are no local points.  Thus $S$ and $T$ are disjoint, so that $T$ is a
finite set of finite places $v$ each of which lies over a prime $p\ge
37$.
To handle this remaining set of places it
suffices to intersect $\cQ_1\cap\cQ_2$ with a suitable 4-plane, as our
next result shows.
\begin{lemma}\label{fin-}
If $T$ is nonempty there is a chain-admissible 4-plane $L'$, 
defined over $K$, such that $L\subset L'$, for which
$\cQ_1\cap\cQ_2\cap L'$ has $K_v$-points for every $v\in T$.
\end{lemma}
It follows of course that $\cQ_1\cap\cQ_2\cap L'$ has points over
every completion of $K$. Naturally, if $T$ were empty the same would
already be true for $\cQ_1\cap\cQ_2\cap L$.

To prove Lemma \ref{fin-} we choose
$K$-points $[\b{e}_0],\ldots,[\b{e}_3]$ spanning $L$, and look for an
additional $K$-point $[\x]=[\b{e}_4]$ such that 
\[L':=<[\b{e}_0],\ldots,[\b{e}_4]>\]
fulfils the necessary conditions. In order to have $L'\subset\cQ_3$ we
will require $[\x]\in\cQ_3\cap(Q_3L)^{\perp}$
We would like the variety $L'$ to be
chain-admissible, and so we will require that $L'\not\in Z_4$.
However $L$ itself is chain-admissible, so that there is at least one
point $P_0$, (which might be defined over $\overline{K}$) for which
\[<L,P_0>\in F_4-Z_4.\]
Such a point $P_0$ will be smooth point of
$\cQ_3\cap(Q_3L)^{\perp}$. It follows that the set of points $P$ for
which $<L,P>$ is a chain-admissible 4-plane, is a nonempty
Zariski-open subset ($U$, say) of $\cQ_3\cap(Q_3L)^{\perp}$. 

To arrange that
$\cQ_1\cap\cQ_2\cap L'$ has a $K_v$-point it will be helpful if $[\x]$
is ``near'' to a $K_v$-point of $\cQ_1\cap\cQ_2$.  We would
therefore like the variety
\[\left(\cQ_3\cap(Q_3L)^{\perp}\right)\cap\left(\cQ_1\cap\cQ_2\right)
=(Q_3L)^{\perp}\cap\cR\]
to contain a $K_v$-point $[\x_v]$ for every $v\in T$. However
$(Q_3L)^{\perp}$ has dimension $n-5\ge 14$, so that
$(Q_3L)^{\perp}\cap\cR$ is the zero locus of a system of three
quadratic forms in at least 15 variables.  There is therefore a 
$K_v$-point whenever $v$ is a finite place above a prime $p\ge 37$, by 
Heath-Brown \cite[Corollary 1]{HBQs}. (The reader should note
that one only needs to know that some bound of the form $p\ge p_0$
suffices. This is a corollary of the famous Ax--Kochen Theorem
\cite{AK}.  However in the present situation one can now provide an
explicit value of $p_0$.  Indeed the work of Schuur \cite{SS} could
also be used here.)  Thus there are points 
$[\x_v]\in\cR\cap(Q_3L)^{\perp}$ for every $v\in T$. Indeed, since $L$
is chain-admissible it is certainly admissible, so that
$\cR\cap(Q_3L)^{\perp}$ must be nonsingular.  If then follows that the
$K_v$-points on $\cR\cap(Q_3L)^{\perp}$ will be Zariski-dense. In
particular we can choose a $K_v$-point $[\x_v]$ in $\cR\cap\M-\cR_1$,
with $\cR_1$ as in Lemma \ref{good4}. Thus, if we set
$L_v=<L,[\x_v]>$, then $L_v$ is admissible so that $\cQ_1\cap\cQ_2\cap
L_v$ is nonsingular.  In particular $[\x_v]$ will be a smooth point of
$\cQ_1\cap\cQ_2\cap L_v$.

Since $[\x_v]\in\cR\cap\M$ it follows
in particular that $[\x_v]\in\cQ_3\cap(Q_3L)^{\perp}$.
We now claim that $[\x_v]$ cannot be a singular point of
$\cQ_3\cap(Q_3L)^{\perp}$. For otherwise the vectors
\[Q_3\x_v,Q_3\b{e}_0,\ldots,Q_3\b{e}_3\]
would be linearly dependent. Since $Q_3$ is nonsingular it would
follow that 
\[\x_v\in<\b{e}_0,\ldots,\b{e}_3>,\]
so that $[\x_v]\in L$.  However $[\x_v]$ was chosen to lie in $\cR$,
so that it would in particular be a $K_v$-point of $\cQ_1\cap\cQ_2\cap L$.
We would therefore have
a contradiction, since $T$ was defined to be the set of places where
$\cQ_1\cap\cQ_2\cap L$ had no local points.

Thus $[\x_v]$ is a smooth point of
$\cQ_3\cap(Q_3L)^{\perp}$.  This variety certainly has at
least one smooth $K$-point, since the quadratic form $Q_3$ was
constructed to split off at least $(n-5)/2>4$ hyperbolic planes over $K$.
Thus we can use weak approximation on $\cQ_3\cap(Q_3L)^{\perp}$ to
produce a sequence of $K$-points $[\x^{(m)}]$, such that 
\[\lim_{m\rightarrow\infty}\x^{(m)}=\x_{v}\]
for every $v\in T$. For large enough $m$ a continuity argument shows that
if $L^{(m)}=<L,[\x^{(m)}]>$ then $L^{(m)}\in F_4-Z_4$.  
We also see
that $\cQ_1\cap\cQ_2\cap L^{(m)}$ approximates $\cQ_1\cap\cQ_2\cap L_v$
in the sense of Lemma \ref{approxV}.  Moreover $[\x_v]$ is a
smooth point of $\cQ_1\cap\cQ_2\cap L_v$, whence Lemma \ref{approxV}
provides $K_v$-points $P^{(m)}$ on $\cQ_1\cap\cQ_2\cap L^{(m)}$ as
soon as $m$ is large enough, for all $v\in T$.  The lemma then follows
on choosing $L'=L^{(m)}$ with a suitably large $m$.
\bigskip

As the final step in our argument we state the
following result.
\begin{lemma}\label{fin}
There is an admissible 7-plane $L''$ containing $L'$, defined over
$K$, and such that
 $\cQ_1\cap\cQ_2\cap L''$ has points over every completion of $K$.
\end{lemma}

Before presenting the proof of the lemma, we show how it suffices for
our theorem. Since $L''\in F_7$ is admissible, the variety $\cQ_1\cap\cQ_2\cap
L''$ is contained in $\cR$, and is a nonsingular intersection of two
quadrics in $\Proj^7$.  By construction it has points over $K_v$ for
every place $v$.  Indeed the subvariety $\cQ_1\cap\cQ_2\cap L$ has a
$K_v$-point $[\b{y}_v]$ for every $v\in S$, and has $K_v$-points for
all $v\not\in T$, while $\cQ_1\cap\cQ_2\cap L'$ has $K_v$-points for
any remaining places $v\in T$. We may therefore apply the
author's result \cite[Theorem 1]{HB8} described in the introduction,
which shows that $\cQ_1\cap\cQ_2\cap L''$ satisfies the Hasse principle
and weak approximation.  This allows us to conclude that
$\cQ_1\cap\cQ_2\cap L''$ has $K$-points arbitrarily close  to
$[\b{y}_v]$ for each $v\in S$.  Our theorem therefore follows.
\bigskip

It remains to establish Lemma \ref{fin}. We have already observed that
either $\cQ_1\cap\cQ_2\cap L$ has points everwhere locally, if $T$ is
empty, or $\cQ_1\cap\cQ_2\cap L'$ does.  The linear spaces $L$ and
$L'$ are chain-admissible, and are defined over $K$. It is therefore
enough to show that if $M$ is any chain-admissible linear space of
dimension $t\in[3,6]$, defined over $K$, then there is a chain
admisssible space $M'\supset M$ of dimension $t+1$, also defined over
$K$. Once this is proved we can use this repeatedly to go from $L$ or
$L'$ to $L''$. 

We recall that $\cQ_3$ contains at least one 7-plane defined over $K$,
whence $M$ will be contained in such a 7-plane, $M^*$, say. We can
choose a basis $\b{e}_1,\ldots\b{e}_n$ for $K^n$ so that $M$ is
spanned by $[\b{e}_1],\ldots,[\b{e}_{t+1}]$ and $M^{*}$ by
$[\b{e}_1],\ldots,[\b{e}_8]$, and such that 
\[Q_3(\sum_1^n
X_i\b{e}_i)=X_1X_9+\ldots+X_8X_{16}+Q(X_{17},\ldots,X_n)\]
for a suitable nonsingular form $Q$.  Then $(Q_3M)^{\perp}$ is spanned by
$[\b{e}_1],\ldots,[\b{e}_8]$ and $[\b{e}_{t+10}],\ldots,[\b{e}_n]$,
and one therefore sees that $[\b{e}_8]$ will be a smooth $K$-point of
$\cQ_3\cap(Q_3M)^{\perp}$. Having shown that there is at least one
smooth $K$-point on $\cQ_3\cap(Q_3M)^{\perp}$ we deduce that the
$K$-points are Zariski-dense, so that there is a $K$-point 
$P\in\cQ_3\cap(Q_3M)^{\perp}-E_t$, in the notation of Lemma \ref{good4}.
We can then complete the proof of Lemma \ref{fin} by taking $M'=<M,P>$.

\bigskip
\bigskip

Mathematical Institute,

Radcliffe Observatory Quarter,

Woodstock Road,

Oxford

OX2 6GG

UK
\bigskip

{\tt rhb@maths.ox.ac.uk}

\end{document}